\documentclass[12pt]{amsart}
\usepackage{amssymb,mathrsfs,upgreek}
\newcommand{\CC}{\mathbb{C}}
\newcommand{\ZZ}{\mathbb{Z}}
\newcommand{\PP}{\mathbb{P}}

\newcommand{\OOO}{\mathscr{O}}

\newcommand{\T}{\operatorname{T}}
\newcommand{\Aut}{\operatorname{Aut}}

\newcommand{\diag}{{\operatorname{diag}}}
\newcommand{\z}{{\operatorname{z}}}

\newcommand{\Syl}{{\operatorname{Syl}}}

\newtheorem{theorem}{Theorem}

\numberwithin{theorem}{section}
\numberwithin{equation}{theorem}

\newtheorem{mtheorem}[theorem]{}

\theoremstyle{definition}
\newtheorem{say}[theorem]{}

\title{A simple proof of the non-rationality 
\\
of a general quartic double solid}
\author{Yuri Prokhorov}
\thanks{Supported by the RSF grant, project No. 14-21-00053 dated 11.08.14}

\address{
Steklov Mathematical Institute
and
\newline\indent
Department 
of Algebra, Moscow State
University
and
\newline\indent
National Research University Higher School of Economics
}
\email{prokhoro@mi.ras.ru}

\begin{document}
\begin{abstract}
The aim of this short note is to give a simple 
proof of the non-rationality of the
double cover of the three-dimensional projective space branched over 
a sufficiently general quartic.
\end{abstract}

\maketitle
\section{Introduction}
Throughout this work the ground field is supposed to be the complex number field $\CC$.

A \textit{quartic double solid} is a projective variety represented as a
the double cover of $\PP^3$ branched along a smooth quartic.
It is known that quartic double solids are unirational but not rational \cite{Beauville1977},
\cite{Tikhomirov1986}, \cite{Voisin1988}, \cite{Clemens1991}.
Moreover, a general quartic double solid is not \textit{stably rational}
\cite{Voisin2015a}.
There are also a lot of results related to rationality problems of \emph{singular} quartic double solids
see e.g. \cite{Artin-Mumford-1972}, \cite{Clemens1983}, \cite{Varley1986}, \cite{Debarre1990},
\cite{Przhiyalkovskij-Cheltsov-Shramov-2015}, \cite{CheltsovPrzyjalkowskiShramov2015b}.

The main result of this note is to give a simple proof of the following

\begin{mtheorem}{\bf Theorem.}\label{theorem-main}
Let $X$ be the quartic double solid branched over the surface 
\begin{equation*}
\label{equation-1}
x_{1}^3x_{2}+ x_{2}^3x_{3}+ x_{3}^3x_{4}+ x_{4}^3x_{1}=0. 
\end{equation*}
Then the intermediate Jacobian $J(X)$ is not a sum of Jacobians of curves.
As a consequence, $X$ is not rational.
\end{mtheorem}

\begin{mtheorem}{\bf Corollary.}\label{corollary-main}
A general quartic double solid is not rational.
\end{mtheorem}

Our proof uses methods of A. Beauville \cite{Beauville2012}, \cite{Beauville2013} and Yu. Zarhin \cite{Zarhin2009}.
The basic idea is to find a sufficiently symmetric variety in the family.
Then the action of the automorphism group provides a good tool to prove 
non-decomposability the intermediate Jacobian into a sum of Jacobians of curves
by using purely \textit{group-theoretic} techniques.
Since the Jacobians and their sums form a
closed subvariety of the moduli space of principally polarized abelian varieties,
this shows that a general quartic double solid is not rational\footnote{Recently
V. Przyjalkowski and C. Shramov used similar method
to prove non-rationality of some double quadrics \cite{PrzyjalkowskiShramov2016}.}.

\section{Preliminaries}
\begin{say}{\bf Notation.}
We use standard group-theoretic notation:
if $G$ is a group, then 
$\z(G)$ denotes its center, $[G,G]$ its derived subgroup,
and $\Syl_p(G)$ its (some) Sylow $p$-subgroup.
By $\zeta_m$ we denote a primitive $m$-th root of unity.
The group generated by elements $\upalpha_1,\upalpha_2,\dots$
is denoted by $\langle \upalpha_1,\upalpha_2,\dots\rangle$.
\end{say}

\begin{say}
Let $X$ be a three-dimensional smooth projective variety with $H^3(X,\OOO_X)=0$
and let $J(X)$ be its intermediate Jacobian regarded as a principally polarized abelian variety
(see \cite{Clemens-Griffiths}). Then $J(X)$ can be written, uniquely up to permutations,
as a direct sum
\begin{equation}
\label{equation-decomposition}
J(X)=A_1 \oplus \dots \oplus A_n,
\end{equation}
where $A_1, \dots, A_p$ are indecomposable
principally polarized abelian varieties (see \cite[Corollary 3.23]{Clemens-Griffiths}). 
This decomposition induces a decomposition of tangent spaces 
\begin{equation}
\label{equation-decomposition-tangent}
\T_{0,J(X)}= \T_{ 0,A_1} \oplus \dots \oplus \T_{0,A_n}
\end{equation}
Now assume that $X$ is acted on by a finite group $G$.
Then $G$ naturally acts on $J(X)$ and $\T_{0,J(X)}$ preserving decompositions \eqref{equation-decomposition}
and \eqref{equation-decomposition-tangent}.
\end{say}

\begin{mtheorem}{\bf Lemma.}\label{lemma-action-80}
Let $C$ be a curve of genus $g\ge 2$ and let $\Gamma\subset \Aut(C)$ be a subgroup
of order $2^k\cdot 5$ whose Sylow $5$-subgroup $\Syl_5(\Gamma)$ is normal
in $\Gamma$. Then the following assertions hold:
\begin{enumerate}
\item\label{lemma-action-80-k=2}
if $k=2$, then $g\ge 3$,

\item\label{lemma-action-80-k=4}
if $k=4$, then $g\ge 6$,

\item\label{lemma-action-80-k=5}\label{faithful-action-on-curve}
if $k=5$, then $g\ge 11$.
\end{enumerate}

\end{mtheorem}
\begin{proof}
Let $C':=C/\Syl_5(\Gamma)$ and $g':=g(C')$. Let $P_1,\dots, P_n\in C'$ be 
all the branch points. By Hurwitz's formula
\[
g+4=5g'+2n.
\]
The group $\Gamma':=\Gamma/ \Syl_5(\Gamma)$ of order $2^k$ faithfully acts on $C'$ and permutes $P_1,\dots, P_n$.
\ref{lemma-action-80-k=2}
Assume that $k=g=2$.
Then $g'=0$, $C'\simeq \PP^1$, and $n=3$.
At least one of the points $P_1,P_2,P_3$, say $P_1$, must be fixed by $\Gamma'$. 
But then $\Gamma'$ must be cyclic (of order $4$) and it cannot leave the set $\{P_1,P_2,P_3\}\subset \PP^1$
invariant. This proves 
\ref{lemma-action-80-k=2}.

\ref{lemma-action-80-k=4}
Assume that $k=4$ and $g\le 5$. Then $g'\le 1$.
If $g'=0$, then $n\in \{3,\, 4\}$ and 
the group $\Gamma'$ of order $16$ acts on $C'\simeq \PP^1$ so that
the set $\{P_1,\dots, P_n\}$ is invariant.
This is impossible. If $g'=1$, then, as above, $\Gamma'$ acts on 
an elliptic curve $C'$ leaving
a non-empty set of $n\le 2$ points is invariant. This is again impossible and the contradiction proves 
\ref{lemma-action-80-k=4}.

\ref{lemma-action-80-k=5}
Finally, let $k=5$ and $g\le 10$. Then $g'\le 2$ and $n\le 7$.
If $g'\le 1$, then we get a contradiction as above.
Let $g'=2$, let $C'\to \PP^1$ the the canonical map, and let 
$\Gamma''\subset \Aut(\PP^1)$ be the image of $\Gamma'$.
Since $\Gamma''$ is a $2$-subgroup in $\Aut(\PP^1)$, it is either cyclic or dihedral.
On the other hand, $\Gamma''$
permutes the branch points $Q_1,\dots,Q_6\in \PP^1$ so that
the stabilizer of each $Q_i$ is a subgroup in $\Gamma''$ of index $\le 4$.
Clearly, this is impossible.
\end{proof}

\section{Symmetric quartic double solid}
\begin{say}
Let $X$ be the quartic double solid as in Theorem \ref{theorem-main}.
Then $X$ is isomorphic to a hypersurface given by
\begin{equation}
\label{equation}
y^2+x_{1}^3x_{2}+ x_{2}^3x_{3}+ x_{3}^3x_{4}+ x_{4}^3x_{1}=0,
\end{equation}
in the weighted projective space 
$\PP:=\PP(1^4,2)$, where 
$x_1,\dots,x_4,y$ are homogeneous coordinates with $\deg x_i=1$, $\deg y=2$.

Let $\upalpha$ be the automorphism of $X$ induced by the diagonal matrix
\[
\diag (1,\, \zeta_{40}^{38},\, \zeta_{40}^{4},\, \zeta_{40}^{26};\, \zeta_{40}^{-1})
\]
and let $\upbeta$ be the cyclic permutation $(1,2,3,4)$ of 
coordinates $x_1,x_2,x_3,x_4$.
Since 
\[
\upbeta\upalpha\upbeta^{-1}= \diag (\zeta_{40}^{26}, 1, \zeta_{40}^{38},\, \zeta_{40}^{4};\, \zeta_{40}^{-1})= 
\diag (1,\, \zeta_{40}^{14},\, \zeta_{40}^{12},\, \zeta_{40}^{18};\, \zeta_{40}^{27})= 
\upalpha^{13},
\]
these automorphisms generate the group 
\[
G= \langle \upalpha,\, \upbeta \mid \upalpha^{40}=\upbeta^4=1,\hspace{3pt} 
\upbeta\upalpha\upbeta^{-1}=\upalpha^{13}\rangle \subset \Aut(X),
\quad G\simeq \ZZ/40 \rtimes\ZZ/ 4.
\]
\end{say}

\begin{mtheorem}{\bf Lemma.}\label{subgroup-index10}
Let $G$ be as above. Then we have
\begin{enumerate}
\item \label{subgroup-index10-1}
$\z(G)=\langle \upalpha^{10}\rangle$ and
$[G,G]=\langle \upalpha^4\rangle$,
\item \label{subgroup-index10-2}
the Sylow $5$-subgroup $\Syl_5(G)$ is normal,
\item \label{subgroup-index10-3}
any subgroup in $G$ of index $10$ contains $\z(G)$.
\end{enumerate}
\end{mtheorem}

\begin{proof}
\ref{subgroup-index10-1} can be proved by direct computations and \ref{subgroup-index10-2} 
is obvious because $\Syl_5(G)\subset \langle\upalpha\rangle$.
To prove \ref{subgroup-index10-3} consider a subgroup 
$G'\subset G$ of index $10$. 
The intersection $G'\cap \langle\upalpha\rangle$ is of index $\le 4$ in $G'$.
Hence $G'\cap \langle\upalpha\rangle$ is a $2$-group of order $\ge 4$ and so 
$\upalpha^{10}\in G'\cap \langle\upalpha\rangle$.
\end{proof}

\begin{mtheorem}{\bf Lemma (cf. \cite[0.1(b)]{Voisin1988}).}\label{Lemma-representation0}
There exists a natural exact sequence
\[
0 \to H^2(X,\Omega_{X}^1 ) \to
H^0(X,-K_X)^\vee {\longrightarrow}\CC \to 0.
\]
\end{mtheorem}
\begin{proof} 
Since $X$ is contained in the smooth locus of $\PP$ and $\OOO_{\PP}(X)=\OOO_{\PP}(4)$, we have 
the following exact sequence
\[
0 \longrightarrow \OOO_X(-4)\longrightarrow \Omega_{\PP}^1|_X \longrightarrow \Omega_{X}^1 \longrightarrow 0,
\]
and so
\[
H^2(X, \Omega_{\PP}^1|_X) \to H^2(X,\Omega_{X}^1 ) \to
H^0(X,\OOO_X(2))^\vee \to H^3(X, \Omega_{\PP}^1|_X)\to 0
\]
The Euler exact sequence for $\PP=\PP(1^4,2)$ has the form
\[
0 \longrightarrow \Omega^1_{\PP} \longrightarrow\OOO_{\PP}(-2)\oplus \OOO_{\PP}(-1)^{\oplus 4} \longrightarrow \OOO_{\PP}
\longrightarrow 0.
\]
Restricting it to $X$ we obtain
$H^2(X, \Omega_{\PP}^1|_X)=0$ and $H^3(X, \Omega_{\PP}^1|_X)=\CC$.
\end{proof}

\begin{mtheorem}{\bf Lemma.}\label{Lemma-representation}
We have the following decomposition of $G$-modules:
\[
\T_{0,J(X)}= V_4\oplus V_4'\oplus V_2,
\]
where $V_4$, $V_4'$ are irreducible faithful $4$-dimensional representations and
$V_2$ is an irreducible $2$-dimensional representation with kernel $\langle \upalpha^{8},\,\upbeta^2\rangle$.
Moreover, $\z(G)$ acts on $V_4$ and $V_4'$ via different characters.
\end{mtheorem}

\begin{proof}
Clearly, $\T_{0,J(X)}\simeq H^0(J(X),\Omega_{J(X)})^\vee \simeq H^2(X,\Omega_{X}^1)$
and by Lemma \ref{Lemma-representation0} we have an injection
$\T_{0,J(X)}\hookrightarrow H^0(X,-K_X)^\vee$.
By the adjunction formula $K_X=(K_{\PP}+X)|_X$ and so
\[
H^0(X,-K_X)\simeq H^0(\PP,\OOO_{\PP}(-K_{\PP}-X)).
\]
Consider the affine open subset $U:=\{x_1x_2x_3x_4\neq 0\}$. 
Then $v=y/x_1^2$ and $z_i=x_i/x_1$, $i=2,3,4$ are affine coordinates in 
$U\subset \{x_1\neq 0\}\simeq \mathbb A^4$.
Let $\upomega$ be the $3$-form
\[
\upomega:= \frac {d z_2\wedge d z_3\wedge d z_4 }{\partial \phi/\partial v} =
\frac {d z_2\wedge d z_3\wedge d z_4 }{2 v},
\]
where $\phi=v^2+z_{2}+ z_{2}^3z_{3}+ z_{3}^3z_{4}+ z_{4}^3$ is the equation of $X$ in $U$.
It is easy to check that for any polynomial $\psi(z_2,z_3,z_4)$ of degree $\le 2$
the element $\psi\cdot \upomega ^{-1}$ extends to a section of $H^0(X, -K_X)$. Thus we have 
\[
 H^0(X, -K_X)\simeq \{ \psi(z_2,z_3,z_4)\cdot \upomega ^{-1} \mid \deg \psi\le 2\}.
\]
It is easy to check that the forms 
\begin{equation}
\label{equation-basis}
\upomega ^{-1}, z_2^2 \upomega ^{-1}, z_3^2 \upomega ^{-1}, z_4^2\upomega ^{-1},
z_2 \upomega ^{-1}, z_2z_3 \upomega ^{-1}, z_3z_4 \upomega ^{-1}, z_4 \upomega ^{-1},
z_3 \upomega ^{-1} , z_2z_4\upomega ^{-1}
\end{equation}
are eigenvectors for $\upalpha$ and $\upbeta$ permutes them.
Moreover, the following subspaces
\begin{itemize}
 \item[]
$W_4= \langle \upomega ^{-1},\hspace{5pt}z_2^2 \upomega ^{-1},\hspace{5pt}z_3^2 \upomega ^{-1},\hspace{5pt}z_4^2\upomega ^{-1}\rangle$,

 \item[]
$W_4'= \langle z_2 \upomega ^{-1},\hspace{5pt}z_2z_3 \upomega ^{-1},\hspace{5pt}z_3z_4 \upomega ^{-1},\hspace{5pt}z_4 \upomega ^{-1}\rangle$,

 \item[]
$W_2= \langle z_3 \upomega ^{-1} ,\hspace{5pt}z_2z_4\upomega ^{-1}\rangle$.

\end{itemize}
are $G$-invariant in $H^0(X, -K_X)$. Moreover, in the basis \eqref{equation-basis}
the element $\upalpha$ acts diagonally:
\begin{equation}
\label{equation-W}
\begin{array}{l}
\upalpha|_{W_4} =\diag (\zeta_{40}^{11},\, \zeta_{40}^{7},\, \zeta_{40}^{19},\, \zeta_{40}^{23}),
\\[4pt]
\upalpha|_{W_4'} =\diag (\zeta_{40}^{9},\, \zeta_{40}^{13},\, \zeta_{40}^{},\, \zeta_{40}^{37}), 
\\[4pt]
\upalpha|_{W_2} =\diag (\zeta_{8}^{3},\, \zeta_{8}^{7}),
\end{array}
\end{equation}
and $\upbeta$ acts on each of these subspaces permuting the eigenspaces of $\upalpha$ cyclically. 
% In particular, $W_4$, $W_4'$, $W_2$ are $G$-invariant. 
Thus $\upalpha^{10}$ acts on $W_4$ (resp., $W_4'$)
via scalar multiplication by $\zeta_4^{3}$ (resp., $\zeta_4$).
Put $V_4:=W_4^\vee$, $V_4':=W_4'^\vee$, $V_2:=W_2^\vee$.
\end{proof}

\section{Proof of Theorem \ref{theorem-main}}

\begin{say}
Assume to the contrary to Theorem \ref{theorem-main} that $J(X)$ is a direct sum of Jacobians of curves, i.e. 
in the unique decomposition \eqref{equation-decomposition} we have
$A_i\simeq J(C_i)$, where $C_i$ is a curve of genus $\ge 1$
and $J(C_i)$ is its Jacobian regarded as a principally polarized abelian variety. 
Let $G_i$ be the stabilizer of $A_i$. There is a natural homomorphism 
$\varsigma_i: G_i \to\Aut (C_i)$.
By the 
Torelli theorem 
$\varsigma_i$ is injective and we have 
\begin{equation}
\label{equation-Aut-C}
\Aut (J(C_i))\simeq 
\begin{cases}
\Aut (C_i)&\text{if $C_i$ is hyperelliptic,}
\\
\Aut (C_i)\times \{\pm 1\}&\text{otherwise.}
\end{cases}
\end{equation}
Let us analyze the action of $G$ on the set $\{A_1,\dots, A_n\}$.
Up to renumbering we may assume that subvarieties $A_1,\dots, A_m$ form one $G$-orbit
(however, the choice of this orbit is not unique in general).
Clearly, $m\in \{1,2,4,5,8,10\}$.
Denote the stabilizer of $A_i$ by $G_i$.
Consider the possibilities for $m$ case by case.
\end{say}

\begin{say}
\label{lemma-invariant-A}
{\bf Case: $m=1$,}
that is, $A_1\subset J(X)$ is a $G$-invariant subvariety. 
Since $\z(G)=\langle \upalpha^{10}\rangle$, the only normal subgroup of order 2 in $G$
is $\langle \upalpha^{20}\rangle$. Hence $G$ cannot be decomposed 
as a direct product of groups of orders $2$ and $80$ (otherwise the order of $\upalpha$ 
would be $20$).
If the action of $G$ on $A_1=J(C_1)$ is faithful, then by \eqref{equation-Aut-C} so is the corresponding
action on $C_1$.  So, the curve $C_1$ of genus $\le 10$ admits 
faithful action of the group $G$ of order $2^{5}\cdot 5$.
This contradicts Lemma \ref{lemma-action-80}\ref{lemma-action-80-k=5}.
Therefore the induced representation on $\T_{0,A_1}$ is not faithful.
By Lemma \ref{Lemma-representation}\ $\T_{0,J(C_1)}= V_2$.
In this case $g(C_1)=2$ and the action of $G$ on $J(C_1)$ induces a faithful action of the group
$\bar G:= G/\langle \upalpha^{8},\,\upbeta^2\rangle$ of order $16$.
Since $C_1$ is hyperelliptic, $\bar G$ is contained in $\Aut(C_1)$.
If $\bar G$ contains the hyperelliptic involution $\tau$, then
$\tau$ generates a
normal subgroup of order 2. In this case $\langle \tau\rangle=[\bar G,\bar G]$
and $\bar G/\langle \tau\rangle$ is an abelian non-cyclic group of order $8$.
But such a group cannot act faithfully on $C_1/\langle \tau\rangle\simeq \PP^1$.
Thus $\bar G$ does not contain the hyperelliptic involution.
In this case the image of the induced action of $\bar G$ on canonical sections
$H^0(C_1,\OOO_{C_1}(K_{C_1}))$ does not contain scalar matrices.
Hence this representation is reducible and so it is trivial on $[\bar G,\bar G]$.
On the other hand, the action of $\Aut(C_1)$ on $H^0(C_1,\OOO_{C_1}(K_{C_1}))$
must be faithful
a contradiction.
\end{say}

From now on we may assume that the decomposition \eqref{equation-decomposition}
contains no $G$-invariant summands.

\begin{say}{\bf Case: $m=5$.}
The subspace $\T_{0,A_1}\oplus\cdots\oplus \T_{0,A_5}\subset \T_{0,J(X)}$ is a $G$-invariant
of dimension $5$ or $10$.
On the other hand, $\T_{0,J(X)}$ contains no invariant subspaces of 
dimension $5$ by Lemma \ref{Lemma-representation}.
Hence, $\T_{0,A_1}\oplus\cdots\oplus \T_{0,A_5}= \T_{0,J(X)}$, $\dim A_i=2$, and
$J(X)= \oplus_{i=1}^{5} A_i$.
The stabilizer $G_i\subset G$ is a Sylow $2$-subgroup that faithfully acts on $C_i$
(because $C_i$ is hyperelliptic, see \eqref{equation-Aut-C}). 
Further, $G_i$ permutes the Weierstrass points $P_1,\dots,P_6\in C_i$.
Hence a subgroup $G_i'\subset G_i$ of index 2 fixes one of them.
In this situation, $G_i'$ must be cyclic.
On the other hand, it is easy to see that $G$ does not contain any elements of order $16$,
a contradiction.
\end{say}

\begin{say}{\bf Case: $m=10$.}
Then $A_1,\dots,A_{10}$ are elliptic curves and
$G_i\subset G$ is a subgroup of index $10$. 
By Lemma \ref{subgroup-index10} each $G_i$ contains $\z(G)$. 
Clearly, $\z(G)$ acts on $\T_{0,A_i}$ via the same character.
Since the subspaces $\T_{0,A_i}$ generate $\T_{0,J(X)}$, the group $\z(G)$ 
acts on $\T_{0,J(X)}$ via scalar multiplication.
This contradicts Lemma \ref{Lemma-representation}.
\end{say}

\begin{say}
\label{orbit=8}
{\bf Case: $m=8$.}
Then $A_1,\dots,A_8$ are elliptic curves and the stabilizer $G_1\subset G$ 
is of order $20$. In particular, the Sylow $5$-subgroup $\Syl_5(G)$ is contained in $G_1$.
Since $\Syl_5(G)$ is normal in $G$, we have $\Syl_5(G)\subset G_i$ for $i=1,\dots,8$.
Since the automorphism group of an elliptic curve contains no order $5$ elements, $\Syl_5(G)$ acts trivially on $A_i$.
Therefore, $\Syl_5(G)$ acts trivially on the $8$-dimensional
$G$-invariant subspace $\T_{0,A_1}\oplus \cdots\oplus \T_{0,A_8}$.
This contradicts Lemma \ref{Lemma-representation}.
\end{say}

\begin{say}\label{orbit=4}{\bf Case: $m=4$.}
The intersection $G_1\cap \langle\upalpha\rangle$ is a subgroup of index $\le 4$
in both $G_1$ and $\langle\upalpha\rangle$. Hence, $G_1\ni \upalpha^{4}$ and so $G_1\supset [G,G]$.
In particular, $G_1$ is normal and $G_1=\cdots=G_4$.
If $\dim A_1=1$, then the element $\upalpha^8$ of order $5$ must
act trivially on elliptic curves $A_i\in 0$, $i=1,\dots,4$.
Therefore, $\upalpha^8$ acts trivially on the $4$-dimensional space $\T_{0,A_1}\oplus\cdots\oplus \T_{0,A_4}$.
This contradicts Lemma \ref{Lemma-representation}.

Thus $\dim A_1=2$. Then $\T_{0,A_1}\oplus\cdots\oplus \T_{0,A_4}=V_4\oplus V_4'$.
An eigenvalue of $\upalpha$ on $\T_{0,A_1}\oplus\cdots\oplus \T_{0,A_4}$ 
must be a primitive $40$-th root of unity (see \eqref{equation-W}).
Hence the group $G_1\cap \langle\upalpha\rangle$ acts faithfully on $\T_{0,A_1}$
and $C_1$ (see \eqref{equation-Aut-C}).
By Lemma \ref{lemma-action-80}\ref{lemma-action-80-k=2} $G_1\cap \langle\upalpha\rangle$ is of order $10$,
i.e. $G_1\cap \langle\upalpha\rangle=\langle \upalpha^4\rangle$ and 
the kernel $N:=\ker (G_1\to \Aut(C_1))$ is of order $4$. Thus $G_1=\langle \upalpha^4\rangle\times N$.
In particular, $G_1$ is abelian.
But then the centralizer $\operatorname{C}(\upalpha^8)$ of $\upalpha^8$ contains $N$ and $\langle\upalpha\rangle$.
Therefore, $\operatorname{C}(\upalpha^8)=G$ and $\upalpha^8\in \z(G)$.
This contradicts Lemma \ref{subgroup-index10}\ref{subgroup-index10-1}.
\end{say}

Thus we have excluded the cases $m=1,4,5,8,10$.
The only remaining possibility is that all the orbits of $G$ on $\{ A_i\}$ 
are of cardinality $2$. 

\begin{say}{\bf Case: $m=2$.}
Then $\dim A_1\le 5$ and $G_1$ is a group of order $80$.
By replacing the orbit $\{A_1,A_2\}$ with another one we may assume that 
$\T_{0,A_1}\oplus \T_{0,A_2}\not \subset V_2$
and so $\T_{0,A_1}\oplus \T_{0,A_2}$ coincides with either $V_4$, $V_4'$, or $V_4\oplus V_4'$. 
In particular, $g(C_1)\ge 2$.
Clearly, 
$G_1\cap \langle\upalpha\rangle$ is of order $40$ or $20$. Hence,
$\upalpha^{2}\in G_1$ and so the group $G_1$ cannot be decomposed as a direct product 
$G_1=\langle \upalpha^{20}\rangle\times H$. By the Torelli theorem $G_1$ faithfully acts on $C_1$.
This contradicts Lemma \ref{lemma-action-80}\ref{lemma-action-80-k=4}.
\end{say}
Proof of Theorem \ref{theorem-main} is now complete.

\begin{proof}[Proof of Corollary \textup{\ref{corollary-main}}]
The Jacobians and their sums form a
closed subvariety of the moduli space of principally polarized abelian varieties.
By Theorem \ref{theorem-main}, in our case, this subvariety does not contain 
the subvariety formed by Jacobians of quartic double solids.
Therefore a general quartic double solid is not rational.
\end{proof}

\subsection*{Acknowledgements.}
The author would like to thank C. Shramov and the referee for useful comments
and  corrections.
\def\cprime{$'$} \def\polhk#1{\setbox0=\hbox{#1}{\ooalign{\hidewidth
  \lower1.5ex\hbox{`}\hidewidth\crcr\unhbox0}}}

% \bibliography{my_ref,specific/inv_groups}
% \bibliographystyle{alpha}
\end{document}